\documentclass[11pt,twoside]{article}
\usepackage{amssymb,amsfonts,amsmath,latexsym,graphicx}
\usepackage[dvips]{epsfig}

\newtheorem{lemma}{Lemma}

\newtheorem{corollary}[lemma]{Corollary}
\newtheorem{theorem}{Theorem}

\newcommand{\beq}{\begin{eqnarray}}
\newcommand{\eeq}{\end{eqnarray}}
\newcommand{\bq}{\begin{equation}}
\newcommand{\eq}{\end{equation}}

\newcommand{\Sum}{\displaystyle \sum}

\newcommand{\Frac}{\displaystyle \frac}
\newcommand{\Inf}{\displaystyle \inf}
\newcommand{\Sup}{\displaystyle \sup}

\newcommand{\email}[1]{{\small E-mail: {\textsf {#1}}}}
\newcommand{\http}[1]{{\small Internet: {\textsf {#1}}}}
\newcommand{\affil}[1]{{\small\sl #1}}

\newcommand{\finprf}{\unskip\null\hfill$\;\square$\vskip 0.3cm}
\def\BbC{{\mathbb C}}

\newcommand{\BbR}{{\mathbb R}}

\newcommand{\un}{{\rm Id}}
\def\square{\sqcup\!\!\!\!\sqcap}
\newcommand{\nrm}[2]{\|{#1}\|_{#2}}
\usepackage{color}
\definecolor{darkred}{rgb}{0.8,0,0}
\definecolor{darkblue}{rgb}{0,0,0.8}

\begin{document}
\title{\sl General results on the eigenvalues of operators with gaps, arising from both ends of the gaps.\\ Application to Dirac operators.}

\author{Jean Dolbeault$^*$, Maria J. Esteban$^*$, Eric S\'er\'e\footnote{Partially supported by European Programs HPRN-CT \# 2002-00277 \& 00282.}\\
\affil{Ceremade (UMR CNRS no. 7534), Universit\'e Paris-Dauphine,}\\
\affil{Place de Lattre de Tassigny, 75775 Paris C\'edex~16, France}\\
\email{dolbeaul, esteban, sere@ceremade.dauphine.fr}\\
\http{http://www.ceremade.dauphine.fr/$\sim$dolbeaul, $\sim$esteban, $\sim$sere}}
\date{\today}\maketitle\thispagestyle{empty}

\abstract{This paper is concerned with  {an extension and reinterpretation } of previous results on the variational characterization of eigenvalues in gaps of the essential spectrum of self-adjoint operators.  {We state} two general abstract results on the existence of eigenvalues in the gap and a continuation principle. Then, these results are applied to Dirac operators in order to characterize simultaneously eigenvalues corresponding to electronic and positronic bound states.}

\section{Introduction}

In \cite{Dolbeault-Esteban-Sere-00B} we proved an abstract result on the
variational characterization of the eigenvalues of operators with gaps in the essential spectrum. Such a result was designed to deal with nonpositive perturbations of a fixed self-adjoint operator with a gap in its essential spectrum but without eigenvalues. In that case, the ``branching'' of the potential ``pulls down'' eigenvalues from the right hand side of the gap. In other words, these eigenvalues emerge from the right end of the gap when the coupling is turned on. Here we address the general case of a perturbation with negative and positive parts, so that eigenvalues can emerge simultaneously from the left and right hand sides of the gap. We thus observe that a simple extension of the general abstract result proved in  \cite{Dolbeault-Esteban-Sere-00B} allows us to treat much more general cases.

 {For a historical overview of the subject, we refer the reader to the introduction of \cite{Dolbeault-Esteban-Sere-00B}, in which an extended review of the literature on eigenvalues in gaps of the essential spectrum is presented.} Some relevant physics papers dealing with this problem are \cite{[Dr-Go], [Ku], [Ta], [Da-De]} (see also the references therein). On the mathematical side we can quote (in chronological order) \cite{[Es-Se], Griesemer-Siedentop-99, Griesemer-Lewis-Siedentop-99, [Dol-Est-Ser], Dolbeault-Esteban-Sere-00B}.

\medskip Let ${\mathcal H}$ be a Hilbert space with scalar product $(\cdot,\cdot)$, and $A:D(A) \subset {\mathcal H} \rightarrow {\mathcal H}$ be a self-adjoint operator. We denote by ${\mathcal H}'$ the dual of ${\mathcal H}$ and by ${\mathcal F} (A)$ the form-domain of $A$. Let ${\mathcal H}_+$, ${\mathcal H}_-$ be two orthogonal Hilbert subspaces of ${\mathcal H}$ such that ${\mathcal H}={\mathcal H}_+ {\displaystyle \oplus} {\mathcal H}_-$. We denote by  $\Lambda^+,\,\Lambda^-$ the projectors on ${\mathcal H}_+$, ${\mathcal H}_-$. We assume the existence of a core $F$ ({\sl i.e.} a subspace of $D(A)$ which is dense for the norm $\nrm\cdot{D(A)}$), such that : \begin{itemize}
\item[ (i)] $F_+ = \Lambda^+ F$ and $F_- = \Lambda^- F$ are two subspaces of ${\mathcal F}(A)$.
\item[(ii$^-$)] $a^-:=\sup_{x_- \in F_-\setminus \{ 0\}} \frac{(x_-, Ax_-)}{\Vert x_- \Vert^2_{_{\mathcal H}}} <+\infty $.
\item[(ii$^+$)]$a^+:=\inf_{x_+ \in F_+\setminus \{ 0\}} \frac{(x_+, Ax_+)}{\Vert x_+ \Vert^2_{_{\mathcal H}}} >-\infty $.
\end{itemize}
We consider the two sequences of min-max and max-min levels $(\lambda_k^+)_{k\geq 1}$ and $(\lambda_k^-)_{k\geq 1}$ defined by
\bq \lambda^+_k := \ \inf_{ \scriptstyle V \ {\rm subspace \ of \ } F_+ \atop \scriptstyle {\rm dim} \ V = k } \ \Sup_{ \scriptstyle x \in ( V \oplus F_- ) \setminus \{ 0 \} } \ \Frac{(x, Ax)}{\|x\|^2_{_{\mathcal H}}}\;,\label{min-max} \eq
\bq \lambda^-_k := \ \sup_{ \scriptstyle V \ {\rm subspace \ of \ } F_- \atop \scriptstyle {\rm dim} \ V = k } \ \inf_{ \scriptstyle x \in ( V \oplus F_+ ) \setminus \{ 0 \} } \ \Frac{(x, Ax)}{\|x\|^2_{_{\mathcal H}}}\;.\label{max-min} \eq
The sequences $(\lambda_k^+)_{k\geq 1}$ and $(\lambda_k^-)_{k\geq 1}$ are respectively nondecreasing and nonincreasing.  {As a consequence of their definitions} we have :
\bq\label{always} \mbox{for all }\; k\geq 1, \quad \lambda^+_k\geq \max\,\{a^-, a^+\} \;\mbox{ and }\; \lambda^-_k\leq \min\,\{a^-, a^+\}\,.
\eq
Let $ b^-\!:=\inf\left\{ \sigma_{{\rm ess}}(A)\cap (a^-, \infty) \right\}$, $b^+\!:=\sup\left\{ \sigma_{{\rm ess}}(A)\cap (-\infty, a^+) \right\}$, and consider the two cases
$${({\rm iii}^-)}\qquad k_0^+:=\min\,\{k\geq 1\,,\;\lambda^+_{k} > a^-\}\,,$$
$${({\rm iii}^+)}\qquad k_0^-:=\min\,\{k\geq 1\,,\;\lambda^-_{k} <a^+\}\,.$$
\begin{theorem}\label{S1}  {If {\rm (i)-(ii$^-$)-(iii$^-$)} hold, for any $k\geq k_0^+$, either $\,\lambda^+_k \,$ is the $(k-k_0+1)$-th eigenvalue of $A$ in the interval $\,(a^-, b^-)\,$ or it is equal to $\,b^-$. If {\rm (i)-(ii$^+$)-(iii$^+$)} hold, for any $k\geq k_0^-$, either $\,\lambda^-_k \, $ is the $(k-k_0+1)$-th eigenvalue of $A$ (in reverse order) in the interval $\,(b^+, a^+)\,$ or it is equal to $\,b^+$.}\end{theorem}
 {Eigenvalues are counted with multiplicity and the order has no meaning if, for instance, $\lambda^+_k=\lambda^+_{k+1}$.}  {The above result does not state anything about the possible eigenvalues of $A$ in the interval $\,[a^+, a^-]$,  {if $a^-\geq a^+$}. We will extensively comment on this in Section 2 and} explain how the abstract result of \cite{Dolbeault-Esteban-Sere-00B} implies Theorem \ref{S1} and a continuation  {result. In Section 3} we will address the particular case in which the operator $\,A\,$ is of the form $\,H_0+V$, where $\,H_0\,$ is the free Dirac operator and $\,V\,$ is an electrostatic scalar potential.

\section{ {Abstract results}}

Theorem 1.1 in \cite{Dolbeault-Esteban-Sere-00B}  {can be stated as follows:}

\smallskip
 {\sl Under the assumptions of Theorem~\ref{S1}, if $\lambda_1^+>a^-$, then all eigenvalues in $(a^-,b^-)$ are given by the min-max levels $\lambda_k^+$ as long as they take their values in $(a^-,b^-)$ (and otherwise, $\lambda_k^+=b^-$).}

\medskip
 {This result dealt with the family of eigenvalues $\,\{\lambda_k^+\}_k\,$ and only in the case $k_0^+=1$. Nothing was said on eigenvalues below $a^+$. The result in \cite{Dolbeault-Esteban-Sere-00B} was already covering all cases corresponding to a Dirac operator with a potential given by a positive Coulomb singularity. Here, by considering the case $\,k_0^+\geq 1$ and by considering the levels $\lambda^-_k$ as well, we extend the method to a framework with interesting physical applications.}

\medskip  {The proof for $\,k_0^+>1\,$ is similar to the proof given in \cite{Dolbeault-Esteban-Sere-00B} and we will not reproduce it here. {\sl A posteriori,\/} passing from $\,k_0^+=1\,$ to $\,k_0^+>1\,$ is not very difficult. Consider indeed a $(k_0^+\!-\!1)$-dimensional} space of  {$F^+$,} $V_{k_0^+\!-\!1}$, such that 
$$\,a^-=\lambda_{k_0^+\!-\!1} \; \leq
\Sup_{ \scriptstyle x \in ( V_{k_0^+\!-\!1} \oplus F_- ) \setminus \{ 0 \} } \ \Frac{(x, Ax)}{\|x\|^2_{_{\mathcal H}}} < \lambda_{k_0^+}\,,$$
and define a new decomposition $\,{\mathcal H}=\tilde{\mathcal H}^+{\displaystyle \oplus} \tilde{\mathcal H}^-\,$ by setting $\,\tilde{\mathcal H}^-={\mathcal H}^- \oplus V_{k_0^+\!-\!1}$.  {Then the first case of Theorem~\ref{S1} is reduced to the result of Theorem 1.1 in \cite{Dolbeault-Esteban-Sere-00B}.}

 {As for the second case}, note that the statement concerning the family $\,\{\lambda_k^-\}_k\,$
follows from that concerning $\,\{\lambda_k^+\}_k\,$ applied to the operator $-A$.  {This completes the sketch of the general ideas for the proof of Theorem~\ref{S1}.\finprf}

\medskip Next, as in \cite{Dolbeault-Esteban-Sere-00B}, we can also consider $\,1$-parameter families of self-adjoint operators of the form $\,A_\tau:=A_0+\tau\,V\,$, $\tau\in [0, \bar\tau]= {{\mathcal I}}\,,$ $\,V\,$ being a bounded scalar potential. In this case, it would be interesting  to prove (iii$^\pm$) for all $\,A_\tau\,$ knowing that $\,A_0\,$ satisfies it and having some spectral  {information on $\,A_\tau$.}

 {More precisely,} let $A_0:D(A_0) \subset {\mathcal H} \rightarrow {\mathcal H}$ be a self-adjoint operator. Let ${\mathcal H}_+$, ${\mathcal H}_-$, $\Lambda^+$ and $\Lambda^-$ be defined like in Section 1. 
Assume further that there is a space $F\,\subset {\cal H}\,$ such that, for all $\,\tau\in {\mathcal I}\,,$ $\,F\,$ is a core for $\,A_\tau$ and the following hypotheses hold:
\begin{itemize}
\item[(j)] $F_+ = \Lambda^+ F$ and $F_- = \Lambda^- F$ are two subspaces of ${\mathcal F}(A_\tau)$.
\item[(jj$^-$)]  {There is $\,a^-\in\BbR$ such that ${\sup_{\tau\in {\mathcal I},\,x_- \in F_-\setminus \{ 0\}} \frac{(x_-, A_\tau\, x_-)}{\Vert x_- \Vert^2_{_{\mathcal H}}} \leq a^-}$.}
\item[(jj$^+$)]  {There is $\,a^+\in\BbR$ such that ${\inf_{\tau\in {\mathcal I},\,x_+ \in F_+\setminus \{ 0\}} \frac{(x_+, A_\tau\, x_+)}{\Vert x_+ \Vert^2_{_{\mathcal H}}} \geq a^+}$.}
\end{itemize}
Let us define the numbers $(\lambda_k^{\tau, +})_{k\geq 1}$ and $(\lambda_k^{\tau, -})_{k\geq 1}$ as in (\ref{min-max})-(\ref{max-min}) by
\begin{eqnarray*}
&&\lambda^{{\tau, +}}_k = \ \inf_{ \scriptstyle V \ {\rm subspace \ of \ } F_+ \atop \scriptstyle {\rm dim} \ V = k } \ \Sup_{ \scriptstyle x \in ( V \oplus F_- ) \setminus \{ 0 \} } \ \Frac{(x, A_\tau)}{\|x\|^2_{_{\mathcal H}}}\,,\label{min-max-nu} \\
&&\lambda^{{\tau, -}}_k = \ \sup_{ \scriptstyle V \ {\rm subspace \ of \ } F_- \atop \scriptstyle {\rm dim} \ V = k } \  {\Inf_{ \scriptstyle x \in ( V \oplus F_+ ) \setminus \{ 0 \} }} \ \Frac{(x, A_\tau)}{\|x\|^2_{_{\mathcal H}}}\,.\label{max-min-nu}\end{eqnarray*}
With the definitions
\begin{eqnarray*}
& {a_1^-:=\inf_{\tau\in {\mathcal I}}\left[\inf\Big( \sigma(A_\tau)\cap(a^-, +\infty)\Big)\right],}&\\
& {a_1^+:=\sup_{\tau\in {\mathcal I}}\left[\sup\Big( \sigma(A_\tau)\cap(-\infty,a^+)\Big)\right],}&\\ \\
&b^{ -}:=\inf\Big(\sigma_{{\rm ess}}(A_0)\cap (a^-, +\infty)\Big),&\\
&b^{ +}:=\sup\Big(\sigma_{{\rm ess}}(A_0)\cap (-\infty, a^+)\Big),&
\end{eqnarray*}
we obtain the following continuation principle.
\begin{theorem}\label{TT5} Under the above assumptions,

\smallskip
 {if for some $\,k_0^+\geq 1$,  $\,\lambda_{k_0^+}^{0,+}>a^-$ and if $\,a_1^->a^-$, for all $\,k\geq k_0^+$, the numbers $\,\lambda^{\tau,+}_k\,$ are either eigenvalues of $\,A_0+\tau\,V\,$ in the interval $\,(a^-, b^{ -})\,$ or $\,\lambda^{\tau,+}_k=b^{ -}\,$.}

{If for some $\,k_0^-\geq 1$, $\, \lambda_{k_0^-}^{0,-}<a^+\,$ and $\,a_1^+<a^+$, for all $\,k\geq k_0^-$, the numbers $\,\lambda^{\tau,-}_k\,$ are either eigenvalues of $\,A_0+\tau\,V\,$ in the interval $\,(b^{ +}, a^+)\,$ or $\,\lambda^{\tau,-}_k=b^{ +}\,$.}
\end{theorem}

 {Exactly as in \cite{Dolbeault-Esteban-Sere-00B}, one can prove this resutl for a class of more general (unbounded) potentials $V$ using a truncation argument and then passing to the limit in the truncation parameter. This applies to the perturbation of free the Dirac operator studied in Section 3 by potentials with Coulomb singularites. We refer the interested reader to \cite{Dolbeault-Esteban-Sere-00B} for more details.}

\medskip\noindent {\it Proof of Theorem \ref{TT5}}. Assumptions (i), (ii$^\pm$) of Theorem \ref{S1} follow from (j), (jj$^\pm$). 
 {Because of the boundedness of $V$, the maps $\,{\mathcal I}\ni\tau\mapsto \lambda_{k_0^\pm}^{\tau, \pm}\;$ are continuous.} 
The sets
$${P_{k_0}^+}:=\{\tau\in {\mathcal I} : \lambda_{k_0}^{\tau,+}\geq a_1^-\}\,,\quad  {P_{k_0}^-}:=\{\tau\in {\mathcal I} : \lambda_{k_0}^{\tau,-}\leq a_1^+\}$$
are thus closed in ${\mathcal I}$, and the sets
$$ {Q_{k_0}^+}:=\{\tau\in {\mathcal I} : \lambda_{k_0}^{\tau,+}> a^-\}\,,\quad  {Q_{k_0}^-}=\{\tau\in {\mathcal I} : \lambda_{k_0}^{\tau,-}< a^+\}$$
are open. Obviously, $P_{k_0}^\pm\subset  {Q_{k_0}^\pm}\,.$ But if $\tau\in  {Q_{k_0}^\pm}$ then $A_\tau$ satisfies (iii$^\pm$), so it follows from Theorem \ref{S1} that
$$\,\lambda_k^{\tau,\pm}\in \sigma(A_\tau)\;, \quad \hbox{for all}\;\,k\geq {k_0}\;,$$
hence, by our assumptions, $\tau\in P_{k_0}^\pm\,.$ As a consequence, $P_{k_0}^\pm= {Q_{k_0}^\pm}$, and the sets  $P_{k_0}^\pm$ are both open and closed in ${\mathcal I}\,.$ But if $\,\lambda_{k_0}^{0,+}>a^-$ (resp. $\,\lambda_{k_0}^{0,-}<a^+$), $ {Q_{k_0}^+}$ (resp. $ {Q_{k_0}^-}$) is nonempty : It contains $0$, so $ {Q_{k_0}^+}$ (resp. $ {Q_{k_0}^-}$) coincides with $\mathcal I$. \finprf

\subsubsection*{ {Example: A Pauli type operator}} For every $\,\nu>0\,$ let us consider the operator

$$ {A_\nu} =\left( \begin{matrix} 1-\Delta-\frac{\nu}{|x|}&0\\ 0&-1+\Delta+\frac{\nu}{|x|}
\end{matrix}\right)\,,$$
on $\,L^2(\BbR^3, \BbC)^2$. This operator is self-adjoint with domain $\,H^2(\BbR^3, \BbC)^2\,$ and form-domain $\,H^1(\BbR^3, \BbC)^2$. An easy analysis shows that for all $\,\nu>0$, $ {A_\nu}$ has two families of eigenvalues:
$$E^-_{\nu, n}=-1+\frac{\nu^2}{4\, {n^2}}\,,\quad E^+_{\nu, n}=1-\frac{\nu^2}{4\, {n^2}}\,,\quad n\geq 1\,,$$
 {and moreover $\,a^\pm_\nu=E^\pm_{\nu, 1}$.}

Furthermore,  {for all $\,k\geq 1$, $\lambda^\pm_{\nu, k}=E^\pm_{\nu, n(k)}$} if and only if 
$\nu\leq\sqrt{\frac{8\, {n^2}}{ {n^2}+1}}\,$,  {$n=n(k)$. Notice indeed that the eigenvalues are degenerate for any $n\geq 2$, so that we have to count the levels with multiplicity and introduce $n:=n(k)$. If $\nu\in\left( \sqrt{\frac{8\, {n^2}}{ {n^2}+1}},\sqrt 8\right)$,}
$\;\lambda^\pm_{\nu,k}=E^\mp_{\nu, 1}\,$ for any $n=n(k)\geq 1$.

\medskip Hence, if $\,\nu\leq 2$,  all the eigenvalues of operator $\, {A_\nu}\,$ are given by the variational procedures defining the numbers $\,\lambda^\pm_k$'s. In the interval $\,\nu\in (2, \sqrt{8})$ some (but not all) of them still satisfy this property. These results are illustrated in Fig.~1 below. 
\begin{figure}[ht]\label{Fig1}\epsfxsize=12.5cm\hbox to\hsize{\epsfbox{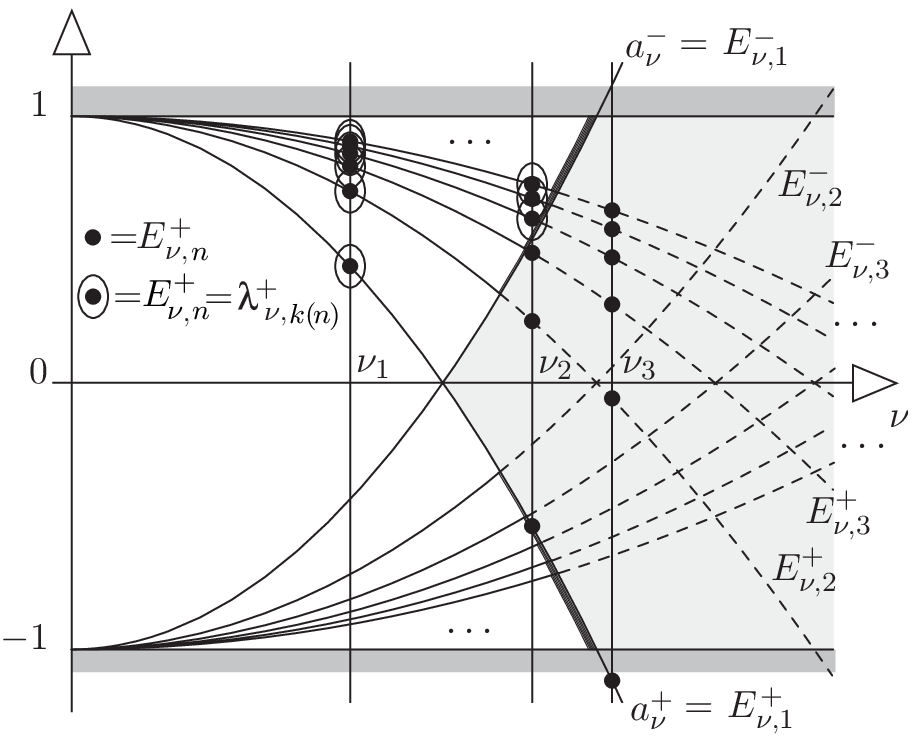}}\caption{\sl Depending on the values of $\nu$, all eigenvalues are achieved by the two families of levels $\lambda^+_{\nu,k}$ and $\lambda^-_{\nu,k}$ (Case $\nu=\nu_1<2$), or only some of them (Case $\nu=\nu_2\in (2,\sqrt 8)$. For $\nu>\sqrt 8$ (Case $\nu=\nu_3$), the gap $(-1,1)$ is contained in $(a^+,a^-)$ and the method does not characterize any eigenvalue in the gap. To clarify the picture, only the family of eigenvalues $E^+_{\nu,n}$ has been represented, but the family $E^-_{\nu,n}$ is easily recovered by symmetry with respect to the horizontal axis.  {To take the multiplicity into account, we denote by $k(n)$ the smallest $k$ for which $\lambda^+_{\nu,k}=E^+_{\nu,n}$.}}\end{figure}

\section{Application to Dirac operators}

Let us consider the free Dirac operator
\[
H_0 := - i\,\sum^3_{k=1} {\bf \alpha}_k\, \partial _k + {\bf \beta} \ , \label{dir-evol}
\]
where we have written it in physical units for which the speed of light, the mass of the electron and Planck's constant are taken equal to $1$. In the Dirac equation, $\alpha_1, \alpha_2, \alpha_3$ and $\beta$ are $4 \times 4$ complex matrices, whose standard form (in $2\times 2$ blocks) is
$$\beta=\left( \begin{matrix} I & 0 \\ 0 & -I \\ \end{matrix} \right) \; ,\;\; \alpha_k=\left( \begin{matrix} 0 &\sigma_k \\ \sigma_k &0 \\ \end{matrix}\right) \qquad (k=1,\, 2,\, 3)\,, $$
and $\sigma_k$, $k=1,\, 2,\, 3$, are the $2\times 2$ complex Pauli matrices: $\sigma_1=\Big( \begin{matrix}\scriptstyle 0 &\scriptstyle 1 \\\scriptstyle 1 &\scriptstyle 0 \\ \end{matrix} \Big)$, $\sigma_2=\Big( \begin{matrix}\scriptstyle 0 &\scriptstyle -i \\\scriptstyle i &\scriptstyle 0 \\ \end{matrix}\Big)$, $\sigma_3=\Big( \begin{matrix}\scriptstyle 1 &\scriptstyle 0\\\scriptstyle 0 &\scriptstyle-1\\ \end{matrix}\Big)$. 
Let $V$ be a scalar potential satisfying
\bq\label{V1} \lim_{|x|\to +\infty}V(x) = 0 \,, \eq
and assume that it is continuous everywhere except at two finite sets of isolated points, $ \,\{x^+_i\}\,$, $\,\{x^-_j\}, \;i=1,\dots I,\; j=1,\dots, J, \,$ where
\bq\label{V2} \begin{array}{ll} 
\lim_{x\to x^+_i} V(x) = +\infty\,, \quad &\lim_{x\to x^+_i} V(x)\, |x-x^+_i|\leq \nu_i\\ &\\
\lim_{x\to x^-_j} V(x)=-\infty\,,\quad &\lim_{x\to x^-_j} V(x)\,|x-x^-_j|\geq -\nu_j
\end{array}\eq
with $\,\nu_i, \nu_j\in (0,1)$ for all $\, i,\, j$. Under the above assumptions, $ H_0 + V$ has a distinguished self-adjoint extension $A$ with domain ${\mathcal D}(A)$ such that
$$ H^1 (\BbR^3, \BbC^4) \subset {\mathcal D} (A) \subset H^{1/2} (\BbR^3, \BbC^4) \,, $$
the essential spectrum of $A$ is the same as that of $H_0\,$:
$$\sigma_{\rm ess} (A) \ = \ (- \infty, - 1] \cup [1, + \infty) \,, $$
 (see \cite{[Th],[Schmin],[Ne],[Kl-Wu]} ). Finally, $\,V\,$ sends $\,{\cal D}(A)\,$ into its dual, since (\ref{V1})-(\ref{V2}) imply that for all $\,\phi\in H^{1/2}(\BbR^3)$, $V\phi\in H^{-1/2}(\BbR^3)$.

\medskip In this section, we shall prove the validity of a variational characterization of the eigenvalues of $\,H_0+V\, $ corresponding to the positive/negative spectral decomposition of the free Dirac operator $\,H_0\,$:
$${\mathcal H}={\mathcal H}^f_+\oplus{\mathcal H}_-^f\,,$$
with $\,{\mathcal H}_\pm^f= \Lambda^f_\pm {\mathcal H}$, where
$$\Lambda^f_+=\chi_{(0, +\infty)}(H_0)= \frac{1}{2}\Big(\un+\frac{H_0}{\sqrt{1-\Delta}}\Big)\,,$$
$$\Lambda^f_-=\chi_{(-\infty, 0)}(H_0)= \frac{1}{2}\Big(\un-\frac{H_0}{\sqrt{1-\Delta}}\Big)\,.$$
This will be done under conditions which are optimal for the potentials satisfying (\ref{V1})-(\ref{V2}) using Theorem \ref{S1} and \ref{TT5}. As already stated in \cite{Dolbeault-Esteban-Sere-00B}, the theorem is optimal in the sense that it covers the optimal range in the case of Coulomb potentials. If we consider the operator $\,A_\tau:=H_0+\tau\,V\,$, $\tau>0$, with $V$ satisfying (\ref{V1})-(\ref{V2}), our variational characterization will provide us with all eigenvalues of $A_\tau$ as long as $\tau$ is not too big.
\begin{theorem}\label{TT67}  {Take a positive integer $\,k_0\,$ and any $k \geq {k_0}\,$ and let $A$ be the self-adjoint extension of $H_0+V$ defined above, where $V$ is a scalar potential satisfying~(\ref{V1})-(\ref{V2}).} 

{If $\,a^-<\lambda^+_{k_0}<1$, then $\,\lambda^+_k \,$ is either an eigenvalue of $\,H_0+V\,$ in the interval $\,(a^-, 1)$, or $\,\lambda^+_k=1$. If additionally $\,V\geq 0$, then $\,a^-=1\,$ and $ \,\lambda^+_k = 1$.}

{If $\, -1<\lambda^-_{k_0}<a^+$, then $\,\lambda^-_k\,$ is either an eigenvalue of $\,H_0+V\,$ in the interval $\,(-1, a^+)\,$ or $ \,\lambda^-_k=-1$. If additionally $\,V\leq 0$, then $\,a^+=-1\,$ and $\, \lambda^-_k = -1$.}\end{theorem}
The sequences $(\lambda_k^+)_{k\geq k_0^\pm}$ and $(\lambda_k^-)_{k\geq 1}$ are respectively nondecreasing and nonincreasing. The spectrum of $A$ contained in \hbox{$\BbR\setminus[a^+, a^-]$~is}
$$(-\infty,-1]\cup\{\lambda_k^\epsilon\;:\; k\geq 1\,,\;\epsilon=\pm\}\cup[1,+\infty)\;,$$
and we do not state anything about the possible eigenvalues in the interval $\,[a^+, a^-]$. As we showed in the previous section, there can be operators for which {\sl all} or {\sl almost all }the eigenvalues lie in the interval $\,[a^+, a^-]$ and thus, they are not given by the variational procedures defining the $\,\lambda^\pm_k$'s.

\medskip
Theorem~\ref{TT67} easily follows from Theorem \ref{S1}. The details of the proof are left to the reader. The continuation argument of Theorem \ref{TT5} applies. Indeed, first one has to truncate the potential at some level $M$, apply Theorem \ref{TT5}, and then pass to the limit when $M$ goes to $+\infty$. It is worth mentionning that by the continuation principle for the Dirac operators $H_0+\tau\, V$, with $V$ satisfying (\ref{V1})-(\ref{V2}), and the definition of $\lambda^{\tau,\pm}_k$,
$$\lim_{\tau \to 0^+}\lambda_k^{\tau, \pm}= \pm 1\,,\quad\mbox{ for all }\; k\geq 1\,.$$
 {Also notice that Talman's decomposition \cite{[Ta],Dolbeault-Esteban-Sere-00B} {\sl i.e.} the decomposition on ``upper'' and ``lower'' two-components spinors, does not apply here, while the spectral decomposition applies.}
\begin{corollary} Under the assumptions of Theorem \ref{TT67}, a sufficient condition for $\,\lambda^+_1\,$ to be in the interval $\,(a^-, 1)\,$ is :
$$\label{V3} c_1-\frac{\nu}{|x|}\leq V\leq c_2\,, \quad c_1,\; c_2\geq 0,\;\, c_1+c_2-1<\sqrt{1-\nu^2}\,.$$ \end{corollary}
\noindent{\it Proof.\/} It is straightforward to check that $\,a^-\leq c_2-1\,$ and
$$\lambda^+_1(V)\geq \lambda^+_1\Big(-\frac{\nu}{|x|}\Big)-c_1 = {\lambda_1}\Big(H_0-\frac{\nu}{|x|}\Big)-c_1=\sqrt{1-\nu^2}-c_1\,.$$
\finprf

Recall that under assumptions (\ref{V1})-(\ref{V2}), for any $k\geq 1$, for the above result to possibly imply that $\,\lambda^\pm_k$ is an eigenvalue we need that
$$\pm a^\mp < 1\quad\mbox{ and }\quad \pm(\lambda_1^\pm-a^\mp)>0\,.$$
To illustrate our results, we end this paper by giving some sufficient conditions for  these inequalities to hold true. Assume that $V$ satisfies \hbox{(\ref{V1})-(\ref{V2})} and can be written as
$$V = -\Sum_{i\in I} V^-_i +\Sum_{j\in J} V^+_j\,,$$
where the $\,V_i^-$'s (resp. the $\,V_j^+$'s) are nonnegative potentials satisfying (\ref{V1})-(\ref{V2}), with a unique singularity at $\,x_i^-$ (resp. at $\,x_j^+$). If
$$\nu_i\; ,\;\nu_j \in {\left[\,0, \,2/\!\left({\scriptstyle \frac{\pi}{2}+\frac{2}{\pi}}\right) \,\right)} \,,\quad \mbox{for all }\, i\in I\,,\; j\in J\,,$$
it follows from \cite{[Tix]} and \cite{[Bur-Ev]} that there are constants $\,\delta_\ell^\pm\in (0, 1)\,$ such that, for all $\, i\in I\,, \; j\in J $,
$$\label{pot-ineq} \delta_i^-\, H_0 - V_i^-\geq 0\;\; \mbox{ in }\; {\cal H}_+\; , \; \;\delta_j^+ H_0 + V_j^-\leq 0\;\; \mbox{ in }\; {\cal H}_- \,,$$
$${\begin{array}{rl}\displaystyle a^- = \kern-12pt\sup_{\begin{array}{c}e\in F_-\\ \nrm e{\mathcal H}=1\end{array}}\kern-12pt \left(H_0+V\right)\leq&\displaystyle  \kern-12pt\sup_{\begin{array}{c}e\in F_-\\ \nrm e{\mathcal H}=1\end{array}}\kern-12pt\Big(H_0+\Sum_{j\in J}V^+_j\Big)\\\displaystyle \leq&\displaystyle \Big(1-\Sum_j \delta^+_j\Big)\kern-12pt\sup_{\begin{array}{c}e\in F_-\\ \nrm e{\mathcal H}=1\end{array}}\kern-12ptH_0 = \Sum_j \delta^+_j-1\,.\end{array}}$$
So, $a^-<1\,$ if 
\bq\label{AA1}\, \Sum_{j\in J} \delta^+_j< 2\,.\eq
Next, let us estimate $\,\lambda^+_1$. For every $\,e_+\in F_+$,
$${\kern-12pt\sup_{\begin{array}{c}e\in [e_+]{\displaystyle \oplus} F_-\\ \nrm e{\mathcal H}=1\end{array}}\kern-12pt(H_0+V) \geq \kern-12pt\sup_{\begin{array}{c}e\in [e_+]\\ \nrm e{\mathcal H}=1\end{array}}\kern-12pt\Big(H_0-\Sum_{i\in I}\,V_i^-\Big)\geq \Big(1-\Sum_i \delta^-_i\Big)\,,}$$
and hence
$$\lambda_1^+\geq 1-\Sum_i \delta^-_i\,.$$
So, finally, $\, \lambda_1^+ > a^- \,\,$ if
\bq\label{AA2} \Sum_{i\in I} \delta^-_i +\Sum_{j\in J} \delta^+_j \,<\,2\,.\eq
Similar computations show that $ \lambda_1^-\!<\!a^+ $ if (\ref{AA2}) holds and \hbox{$a^+\!>\!-1$} if  {additionally}
\bq\label{AA3}\, \Sum_{i\in I} \delta^-_i< 2\,.\eq

\medskip Conditions (\ref{AA1}), (\ref{AA2}) and (\ref{AA3}) are very restrictive. If the interdistances between the singularity points $\,x_i^-\,$ and $\,x_j^+$ are taken into account and made large enough, these conditions can certainly be radically weakened when these interdistances become large.


\end{document}